\documentclass[12pt]{amsart}

\usepackage{amsmath,amsthm,amsfonts,amssymb} 

\usepackage[all]{xy}

\usepackage{color}

\allowdisplaybreaks

\numberwithin{equation}{section}

\newtheorem{thm}[equation]{Theorem}

\newtheorem{example}[equation]{Example}

\newtheorem{remark}[equation]{Remark}

\newtheorem{remarks}[equation]{Remarks}

\newtheorem{question}[equation]{Question}

\DeclareMathOperator{\gr}{gr}

\DeclareMathOperator{\coh}{H}

\DeclareMathOperator{\Ext}{Ext}

\DeclareMathOperator{\Hom}{Hom}

\DeclareMathOperator{\ima}{Im}

\DeclareMathOperator{\Ker}{Ker}

\DeclareMathOperator{\Span}{Span}

\newcommand{\N}{{\mathbb N}}

\newcommand{\ot}{\otimes}

\newcommand{\res}{{\rm res}}

\begin{document}

\title[Finite generation of cohomology]{Finite generation of the cohomology\\
of some skew group algebras}



\author{Van C.\ Nguyen}

\address{Department of Mathematics\\Hood College\\Frederick, MD 21701}

\email{nguyen@hood.edu}

\author{Sarah Witherspoon}

\address{Department of Mathematics\\Texas A\&M University\\College Station, TX 77843}

\email{sjw@math.tamu.edu}

\date{October 21, 2017}

\thanks{This material is based upon work done while the first author was a Texas A\&M graduate student. It was supported by the National Science Foundation under grant No.\ 0932078000 while both authors were in residence at the Mathematical Sciences Research Institute (MSRI) in Berkeley, California, during the semester of Spring 2013. Both authors were also supported by NSF grant DMS-1101399.}

\maketitle


\begin{abstract}

We prove that some skew group algebras have Noetherian cohomology rings, a property inherited from their component parts. The proof is an adaptation of Evens' proof of finite generation of group cohomology. We apply the result to a series of examples of finite dimensional Hopf algebras in positive characteristic.

\end{abstract}


\section{Introduction}

The cohomology ring of a Hopf algebra encodes potentially useful information about its structure and representations. It is always graded commutative (see, for example, \cite{Suarez-Alvarez}). For many classes of finite dimensional Hopf algebras, it is also known to be finitely generated: for example, cocommutative Hopf algebras (Friedlander and Suslin \cite{Friedlander-Suslin}), small quantum groups (Ginzburg and Kumar \cite{Ginzburg-Kumar}), and small quantum function algebras (Gordon \cite{Gordon}). Etingof and Ostrik \cite{Etingof-Ostrik} conjectured that it is always finitely generated, as a special case of a conjecture about finite tensor categories. Snashall and Solberg \cite{Snashall-Solberg} made an analogous conjecture for Hochschild cohomology, of finite dimensional algebras, that was seen to be false when Xu \cite{Xu} constructed a counterexample. In contrast, there is neither a counterexample nor a proof of  the  Hopf algebra conjecture. Each finite generation result so far has used, in crucial ways, known structure of a particular class of Hopf algebras. Further progress will require new ideas.

In this article, we present one technique for handling some types of algebras inductively. Many (Hopf) algebras of interest are skew group algebras (that is, smash products with group algebras). Under some conditions on a skew group algebra, we show that its cohomology is Noetherian if the same is true of the underlying algebra on which the group acts.

Specifically, if $A$ is a finite dimensional augmented algebra over a field $k$, with an action of a finite group $G$ by automorphisms, there is a  spectral sequence relating the cohomology of the smash product $A\# kG$ (definition in Section~\ref{defns}) as an augmented algebra to that of each  of $A$ and  $G$. (It is essentially the Lyndon-Hochschild-Serre spectral sequence.) This allows us to use the  framework of Evens' classic proof of finite generation of group cohomology \cite{Evens} to prove that the cohomology rings of some smash products are Noetherian (Theorem~\ref{main}). In order to do this, we need a particularly nice set of permanent cycles in the cohomology  of $A$. In the finite group case, these cycles exist due to an application of  Evens' norm map. In our setting, there may be no such norm map, and we instead hypothesize existence of these permanent cycles.

We focus on a class of examples (in Section~\ref{sec:pointed examples}) found by Cibils, Lauve, and the second author~\cite{CLW} that satisfy our hypotheses. We prove finite generation of the cohomology of these noncommutative, noncocommutative Hopf algebras in positive characteristic. While our main theorem is tailored to suit these examples, we state and prove it in the abstract setting, in order to add one more tool to the collection of techniques available for proving finite generation. Our restrictive hypotheses serve to highlight the difficulty in adapting methods designed for the finite group setting, where serendipity reigns.

We thank D.\ Benson and P.\ Symonds for very insightful conversations and suggestions. We thank \O.\ Solberg for computing the cohomology of some of the Nichols algebras in Section~\ref{anick}; these computations led us to our general result on this series of Nichols algebras and corresponding Hopf algebras.


\section{Definitions and notation} 

\label{defns}

Throughout this article, let $k$ be a field. All algebras will be associative algebras over $k$, and all modules will be left modules, finite dimensional over $k$. Let $\ot = \ot_k$.

Let $G$ be a finite group acting on a finite dimensional augmented $k$-algebra $A$ by automorphisms. Let $A\# kG$ be the resulting {\em smash product} (or {\em skew group algebra}), that is $A\ot kG$ as a vector space, with multiplication $(a\ot g)(b\ot h) = a ( ^g \! b ) \ot gh$, for all $a,b\in A$ and $g,h\in G$. (For simplicity, we will drop tensor symbols in this notation from now on.) We assume the action of $G$ preserves the augmentation of $A$, so that $A \# kG$ is also augmented with augmentation map $\varepsilon_{A \# kG}: A \# kG \rightarrow k$ defined by $\varepsilon_{A \# kG}(ag) = \varepsilon_A(a)$, for all $a \in A$, $g \in G$.

We use the symbol $k$ also to denote the one-dimensional $A$-module (respectively, $A \# kG$-module) on which $A$ (respectively, $A \# kG$) acts via its augmentation. Let 
$$\coh^*(A,k):= \Ext^*_A(k,k) \ \ \ \mbox{ and } \ \ \  \coh^*(A \# kG,k) := \Ext^*_{A \# kG} (k,k).$$ 
Both are algebras under Yoneda composition. The embedding of $A$ into $A\# kG$ as a subalgebra induces a {\em restriction map} 
$$\res_{A\# kG , A} : \coh^*(A \# kG,k)\rightarrow \coh^*(A,k)$$ on cohomology. There is an action of $G$ on $\coh^*(A,k)$ that may be defined for example via the diagonal action of $G$ on the components of the bar resolution for $A$. There is a similar action of $G$ on $\coh^*(A \# kG,k)$ that is trivial since it comes from inner automorphisms on $A \# kG$.


\section{Finite generation of cohomology}

In this section, we prove our main theorem that under certain hypotheses, the cohomology ring $\coh^*(A\# kG,k)$ of $A\# kG$ is Noetherian:

\begin{thm}
\label{main}

Let $G$ be a finite group acting on a finite dimensional augmented algebra $A$, preserving the augmentation map. Assume that $\ima (\res_{A\# kG,A})$ contains a polynomial subalgebra over which $\coh^*(A,k)$ is Noetherian and free as a module, with a free basis whose $k$-linear span is a $kG$-submodule of $\coh^*(A,k)$. Then $\coh^*(A\# kG,k)$ is Noetherian.

\end{thm}

\begin{remarks}{\em 
 \item \hspace{0.5cm}(a) The hypothesis that $\ima (\res_{A\# kG,A})$ contains a polynomial subalgebra over which $\coh^*(A,k)$ is Noetherian, together with the left module version of \cite[Corollary 1.5]{GW}, implies that $\coh^*(A,k)$ is (left) Noetherian. 

 \item \hspace{0.5cm}(b) We did not specify the characteristic of the base field $k$ in the theorem. If the characteristic of $k$ does not divide the order of $G$, then $kG$ is semisimple and its cohomology is trivial except in the degree $0$. In this case, $\coh^*(A\# kG,k) \cong \coh^*(A,k)^G$, the invariant ring under the action of $G$. Here, one can use invariant ring theory in the noncommutative setting to show that the conclusion of the theorem  holds. (See, for example, \cite[Corollary~4.3.5]{Mo}.) For the proof of Theorem~\ref{main}, we assume the characteristic of $k$ divides the order of $G$.
}
\end{remarks}

\begin{proof}

We use the Lyndon-Hochschild-Serre spectral sequence (see, for example, \cite[Chapter VI]{Barnes} in a very general setting):
\[E_2^{p,q}=E_2^{p,q}(k)=\coh^p(G,\coh^q(A,k)) \Longrightarrow \coh^{p+q}(A\# kG,k).\]
Let $E_r(k)$ denote the resulting $r$th page, and note that for each $q$, $\coh^q(A,k)$ is a finite dimensional $k$-vector space.

Note that $E_\infty^{0,*}$ is a submodule of $E_2^{0,*}$, since no $d_r: E_r^{p,q} \rightarrow E_r^{p+r,q-r+1}$ ends on the vertical edge. It follows that the restriction map $\coh^*(A\# kG,k) \rightarrow E_2^{0,*}(k)$ is part of the following commuting diagram: 

\begin{displaymath}
  \xymatrix{ 
	   \coh^*(A\# kG,k) \ar[d]       \ar[rrr]  ^{\res_{A\# kG, A}}  &&& \coh^0(G,\coh^*(A,k))=\coh^*(A,k)^G \ar@{=}[d] \\            
	   E_\infty^{0,*}(k) \ar[rrr]^{\hookrightarrow}  &&&  E_2^{0,*}(k) }
\end{displaymath}

We can identify $E_\infty^{0,*}$ with the image of the restriction map in $E_2^{0,*}$.

Let $T=k[\chi_1,\dots,\chi_m]$ denote the polynomial subalgebra of $\ima (\res_{A\# kG, A})$ hypothesized in the statement of the theorem. The action of $G$ on $\coh^*(A,k)$ restricts to the trivial action on $T$ since it is a subalgebra of $\ima (\res_{A\# kG,A})$. Therefore,  by the Universal Coefficients Theorem, $\coh^*(G,T)\cong \coh^*(G,k)\ot T$, an isomorphism of graded algebras.

Let $S:=\coh^{*}(G,k)=E_2^{{*},0}(k)$. Let $R$ be the subring of $E_2(k)$ generated by $S$ and $T$. By the above observations,  $R\cong S[\chi_1,\ldots,\chi_m]$, a polynomial ring over $S$ in $m$ indeterminates (that we also denote by $\chi_1,\ldots,\chi_m$ for convenience). Since $d_2$ vanishes on the horizontal edge,  $R \subseteq \Ker(d_2)$. So $R$ projects onto a subring of $E_3(k) = \coh(E_2(k),d_2)$. Similarly, $R$ projects onto a subring of $E_r(k)$ for every $r>0$ including $\infty$. Therefore, we may consider $E_r(k)$ to be a module over $R$, for every $r>0$ including $\infty$.

\quad

\noindent

{\bf Claim 1: $E_2(k)$ is a Noetherian module over $R$}.

\noindent

{\it Proof of Claim 1.} By hypothesis, there are (homogeneous) elements $\rho_1,\ldots,\rho_t\in\coh^*(A,k)$ that form a free basis of $\coh^*(A,k)$ as a $T$-module, and for which
$$V:=\Span_k\{\rho_1,\ldots,\rho_t\}$$ is a $kG$-submodule of $\coh^*(A,k)$. 
Let 
\[ L:= \coh^*(G, V).\]
Note that $L$ contains a copy of $S=\coh^*(G,k)$ as $V$ must include an element in degree 0, that is in $\coh^0(A,k)\cong k$, which has trivial $G$-action. By hypothesis, $\coh^*(A,k) = k[\chi_1,\ldots,\chi_m]\cdot V$, and so 
$$ E_2(k) = \coh^*(G, \ k[\chi_1,\ldots,\chi_m]\cdot V).$$
Further, $k[\chi_1,\ldots,\chi_m]$ has trivial $G$-action and the  module $\coh^*(A,k)$ for this polynomial ring is free with free basis $\rho_1,\ldots,\rho_t$. It follows that, as a $kG$-module, 
\[  k[\chi_1,\ldots,\chi_m]\cdot V\cong \bigoplus_{i_1,\ldots, i_m\geq 0}  \chi_1^{i_1}\cdots\chi_m^{i_m}\cdot V \cong \bigoplus_{i_1,\ldots,i_m \geq 0}   V, \]
a direct sum of copies of the same $kG$-module, $V$. Therefore by the Universal Coefficients Theorem, $E_2(k)$  is the image of 
$$
   \coh^0(G, k[\chi_1,\ldots,\chi_m])\ot \coh^*(G, V) \cong k[\chi_1,\ldots,\chi_m]\ot L,$$ 
under cup product. We thus identify $E_2^{*,*}(k)$ with $S[\chi_1, \ldots, \chi_m] \ot_S L$.

Since $G$ is a finite group and $V$ is a finite dimensional vector space over $k$, $L=\coh^*(G,V)$ is Noetherian over $S=\coh^*(G,k)$ \cite{Evens}. By the Hilbert Basis Theorem for graded commutative rings (see, for example, \cite[Theorem 2.6]{GW}), $S[\chi_1,\ldots,\chi_m] \ot_S L$ is Noetherian over $R=S[\chi_1,\ldots,\chi_m]$. Therefore, $E_2^{*,*}(k)$ is Noetherian over $R$. We have proven Claim 1.

\quad

\noindent

{\bf Claim 2: The spectral sequence stops, i.e., $E_r=E_\infty$ for some $r<\infty$}.

\noindent

{\it Proof of Claim 2.} Let $Z_i$ be the space of $i$-cocycles and $B_i$ be the space of $i$-coboundaries in $E_i=E_i(k)$. Recall that $E_1=Z_1$ and $E_2=Z_2/B_2$. Consider the ``pull back'' $B_r$ in $E_2$ of $d_r(E_r)$ as follows:

Each element of $E_2$ on which $d_2$ vanishes determines an element of $E_3$. Suppose $d_3$ vanishes on that element, so that it in turn determines an element of $E_4$. Continue placing such restrictions until we determine an element of $E_r$, and suppose that element is in the image of $d_r$. We define:
\[B_r:=\left\{\tau \in E_2: \tau \in \Ker(d_i) \text{, for } 2 \leq i\leq r-1 \text{, and } \tau \in \ima(d_r)\right\}. \]
Note that $B_r$ is an $R$-submodule of $E_2$ since $d_j$ is a derivation for all $j$, $2 \leq j \leq r$, and the image in each $E_j$ of $R$ consists of universal cycles. Moreover, $B_r \subseteq B_{r+1}$ so we obtain an ascending chain of $R$-submodules of $E_2$: 
\[ 0=B_1 \subseteq B_2 \subseteq \cdots \]
Since $E_2$ is Noetherian over $R$ by Claim 1, this chain must stabilize by the ascending chain condition. Thus there exists some $r_0$ finite such that $B_{r_0}=B_{r_0+1}=B_{r_0+2}= \ldots \ $, and so $d_r=0$ for all $r >r_0$. This implies $E_r=E_\infty$ for all $r>r_0$. We have proven Claim 2.

\quad

We can put this together to finish the proof of the theorem: Each $Z_r, B_r$ is a submodule of $E_2$ over $R=S[\chi_1,\ldots,\chi_m]$. Thus, each $E_r$, which is a submodule of a quotient module of $E_{r-1}$, is Noetherian over $R$ by Claim 1 and induction on $r$. By Claim 2, $E_\infty$ is Noetherian over $R$, and so by \cite[Corollary~1.5]{GW} it is a Noetherian ring.

Now, $\coh^*(A\# kG,k)$ has a filtration whose filtered quotients are 
\[E_\infty^{p,q}(k) \cong \frac{F^p\coh^{p+q}(A\# kG,k)}{F^{p+1}\coh^{p+q}(A\# kG,k)}.\]
Suppose that $\coh^*(A\# kG,k)$ is not Noetherian and let $T_1 \subseteq T_2 \subseteq \cdots \subseteq \coh^*(A\# kG,k)$ be an infinite ascending chain of ideals of $\coh^*(A\# kG,k)$. Let 
\[F^pT_i:=T_i \cap F^p\coh^*(A\# kG,k)\] 
and 
\[U_i:= \bigoplus_{p\geq0} F^pT_i/F^{p+1}T_i \subseteq E_\infty(k).\]
If $x \in T_{i+1} \setminus T_i$, then for some $p$, $x \in F^pT_{i+1}$ but $x \notin F^pT_i$ and $x \notin F^{p+1}T_{i+1}$, so $x+F^{p+1}T_{i+1}$ is not in the image of the inclusion 
\[F^pT_i/F^{p+1}T_i \hookrightarrow F^pT_{i+1}/F^{p+1}T_{i+1}\]
that is, $x \in U_{i+1} \setminus U_i$. So $U_{i+1}$ properly contains $U_i$, for all $i$. Therefore, we have an infinite ascending chain of ideals of $E_\infty(k)$:
\[U_1 \subsetneqq U_2 \subsetneqq \cdots \]
This contradicts the result that $E_\infty(k)$ is Noetherian. Hence, $\coh^*(A\# kG,k)$ is Noetherian. \end{proof}

\begin{remark}
{\em Theorem~\ref{main} parallels the main step in Evens' proof of finite generation of group cohomology: Let $H$ be a finite $p$-group (where $k$ has characteristic $p$), $A = kZ$ is the group algebra of a central subgroup $Z$ of $H$ of order $p$, and $G=H/Z$. (In case $Z$ is complemented in $H$, we obtain $kH\cong A\#  kG$, whereas more generally, $kH$ is a crossed product of $A$ with $G$.) In this case, Evens' norm map is applied to show that $\ima (\res_{kH, kZ})$ contains a polynomial subalgebra $k[\zeta]$ (in one indeterminate). One observes that $\coh^*(kZ, k)$ is a free module over $k[\zeta]$, and that the $k$-linear span of any free basis is a $kG$-submodule. This special case is somewhat simpler than our more general context as it uses a polynomial ring in one indeterminate. }
\end{remark}

We are particularly interested in those actions of finite groups $G$ on algebras $A$ for which $A\# kG$ is a Hopf algebra. We turn to a class of such examples in the remainder of the paper.


\section{Examples: Nichols algebras in positive characteristic}

\label{anick}

In this section, we first recall the  Nichols algebras $A$ from \cite[Corollary~3.14]{CLW} and the corresponding Hopf algebras $A\# kG$ from the same paper. We will prove that these Hopf algebras have finitely generated cohomology. This will follow from Theorem~\ref{main} and explicit calculation using Anick's resolution \cite{Anick}. In this section we explain these calculations for $A$, and in the next we complete the proof of finite generation of cohomology of $A\# kG$. The results of this section were anticipated by Solberg~\cite{OS} as a consequence of  computer calculations (for small $p$) that gave the graded vector space structure and generators of cohomology.

In the remainder of the paper, $k$ will be a field of characteristic $p>2$. (The case $p=2$ is included in \cite{CLW}, but is different, and we will not consider that case here.) Let $A$ be the augmented $k$-algebra generated by $a$, $b$, with relations
\[  a^p=0, \ \  \ b^p = 0, \ \ \ ba = ab+ \frac{1}{2} a^2, \]
and augmentation $\varepsilon: A\rightarrow k$ given by $\varepsilon(a)=\varepsilon(b)=0$.
Let $G$ be a cyclic group of order $p$ with generator $g$, acting on $A$ on the right by 
\[ g(a) = a , \ \ \  g(b) = a+b. \footnote{To apply results in \cite{NWW}, we changed from the left $G$-module structure: $g(a) = a, \;  g(b) = b-a$ to this right $G$-module structure.}\]
Then $A\# kG$ is a  Hopf algebra with comultiplication given by 
\[ \Delta(g) = g\ot g, \ \ \ \Delta(a)=a\ot 1 + g\ot a, \ \ \   \Delta(b) = b\ot 1 + g\ot b. \]

It is useful to consider $A$ as a quotient of a larger algebra. Let 
\begin{equation}\label{eqn:B}
   B := k\langle a,b\rangle /(ba - ab - \frac{1}{2}a^2 ), 
\end{equation}
so that $A \cong B/(a^p, b^p)$. We will show that $B$ is a PBW algebra in the sense of \cite{BGV} or \cite[Section 2]{Shroff}, although we will not need this fact for our cohomology calculations.

Choose the lexicographic order on $\N^2$ for which $(0,1) < (1,0)$, and assign $\deg(a) = (0,1), \ \deg(b) = (1,0)$. Then $ba-ab-\frac{1}{2}a^2$ is a Gr\"obner basis for the ideal of the free algebra $k\langle a,b\rangle$ that it generates. It follows that $\{a^ib^j\mid i,j\geq 0\}$ is a vector space basis of $B$. The relation $ab = ba - \frac{1}{2}a^2$ satisfies the required condition in the definition of a PBW algebra since $\deg(a^2) < \deg(ab)$, so $B$ is a PBW algebra. Moreover, $B$ is a Koszul algebra by a theorem of Priddy~\cite[Theorem~5.3]{Priddy}.

Applying \cite[(3.9)]{CLW}, one finds that the elements $a^p, b^p$ are in the center of $B$. We may thus  apply a theorem of Shroff, \cite[Theorem~4.3]{Shroff}, to the Nichols algebra $A$ to conclude that the cohomology ring $\coh^*(A,k)$ of $A$ is Noetherian.

We will need some details about this cohomology of $A$ for the next section. For this, we will construct Anick's resolution \cite{Anick} for $A$, and show that it is minimal. We use the combinatorial description of the resolution given by Cojocaru and Ufnarovski~\cite{CU}, however we  index differently, and use left modules instead of right. This is a free resolution of the trivial $A$-module $k$, of the form 
$$  \cdots \longrightarrow A\ot kC_2 \longrightarrow A\ot kC_1 \longrightarrow A\ot kC_0 \longrightarrow k \rightarrow 0, $$
for (finite) sets $C_n$, where $kC_n$ denotes the vector space with basis $C_n$. Let $C_{0}:= \{1\}$ and $C_1:= \{a,b\}$. Then $C_2:=\{ a^p, \ b^p, \ ba\}$ is the set of ``tips'' or ``obstructions.'' To define $C_n$ in general, consider the graph 

$$
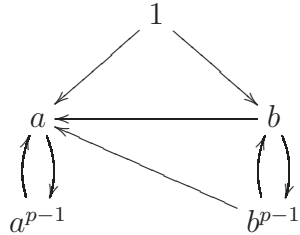
\xymatrix{ &  1 \ar[dl]  \ar[dr] &  \\ 
a\ar@/^/[d]  & & b\ar[ll]\ar@/^/[d] \\
a^{p-1}\ar@/^/[u] & & b^{p-1}\ar[llu]\ar@/^/[u] } $$

\quad

\noindent
The elements of $C_n$ correspond to paths of length $n$ that start at $1$. We label such paths with the product of all elements through which the path passes (including the endpoint). In this way we obtain 

\begin{eqnarray*}
  C_3 & = & \{ a^{p+1}, \ b^{p+1}, b^pa, \ ba^p  \} , \\
  C_4 & = & \{ a^{2p}, \ b^{2p}, b^{p+1}a, \  b^pa^p, \ b a^{p+1} \} ,  
\end{eqnarray*}
and in general
\begin{eqnarray*}
  C_{2m-1} & = & \{ b^{kp} a^{(m-1-k)p + 1}, \ b^{kp+1} a^{(m-1-k)p} \mid k=0,1,\ldots, m-1\} , \\
  C_{2m} & = & \{ b^{mp}, \ b^{kp} a^{(m-k)p}, \ b^{kp+1} a^{(m-1-k)p +1} \mid k=0,1,\ldots, m-1\} .
\end{eqnarray*}

For qualitative understanding of the differentials, give each of the generators $a,b$ of $A$ the degree~$1$. We claim that the differentials preserve degree, where the graded module structure of a tensor product $A\ot kC_i$ is given by $\deg (a\ot x) = \deg (a) + \deg(x)$ if $a,x$ are homogeneous. This claim results from the recursive definition of the differential $d$ in each homological degree: By construction, $d$ applied to elements of $A\ot kC_1$ is multiplication, and to $A\ot kC_2$ takes each tip to the Gr\"obner basis element to which it corresponds, suitably expressed as an element of $A\ot kC_1$. The remaining differentials are defined iteratively, via splitting maps in each homological degree that are also defined iteratively. Since the relations are homogeneous and differentials in low homological degrees preserve degrees of elements, the splitting maps and differentials in higher degrees may be chosen to have the same property.

Now note that $C_{2m-1}$ consists of elements of degree $(m-1)p+1$, and $C_{2m}$ consists of elements of degrees $m p$ and $(m-1)p+2$. Therefore elements of $C_n$ and of $C_{n-1}$ never have the same degree. As a consequence the differential takes elements of $C_n$ to elements of $A_+\ot C_{n-1}$ where $A_+$ denotes all elements of $A$ of positive degree (and these are in the kernel of the augmentation map $\varepsilon$). When applying the functor $\Hom_A( - , k)$ then, the induced differentials all become 0. Therefore in this case, Anick's resolution is minimal, and for each $n$, the dimension of $\coh^n(A,k)$ is $n+1$.


\section{Examples: Pointed Hopf algebras in positive characteristic}
\label{sec:pointed examples}

We wish to apply Theorem \ref{main} to the Hopf algebras $A\# kG$ introduced in the previous section. In order to do this, we next give some of the details from Shroff~\cite[Section 4]{Shroff} as they apply to these examples in particular. Recall the PBW algebra $B$ defined in (\ref{eqn:B}). Let $\xi_a,\xi_b: B\ot B\rightarrow k$ be the $k$-linear functions given by
\[ \xi_a(r\ot s) = \gamma_a, \ \ \  \xi_b(r\ot s) = \gamma_b, \]
where $\gamma_a$ (respectively, $\gamma_b$) is the scalar coefficient of $a^p$ (respectively, $b^p$) in the product $rs$ in $B$. (These functions are denoted $\widetilde{\zeta}_1, \widetilde{\zeta}_2$ in \cite{Shroff}.) Extending to left $B$-module homomorphisms in $\Hom_B(B^{\ot 3},k)$ under the isomorphism $\Hom_B(B^{\ot 3},k)\cong \Hom_k(B^{\ot 2},k)$, the functions $\xi_a,\xi_b$ are coboundaries on the bar resolution of $B$, as shown in \cite{Shroff}, and they factor through $A \cong B/(a^p,b^p)$. The resulting functions (which we will also denote $\xi_a,\xi_b$ by abuse of notation) are no longer coboundaries. They represent nonzero elements in the cohomology of $A$, corresponding to permanent cycles in the May spectral sequence  for $A$ as a  filtered algebra (see \cite[Theorem 3]{May} or \cite[5.4.1]{Weibel}). On page $E_1$ of this spectral sequence, their counterparts generate a polynomial ring over which $E_1$ is finitely generated (by the elements $ 1 , \eta_a, \eta_b, \eta_a\eta_b$, where $\eta_a,\eta_b$ have cohomological degree~$1$, functions dual to $a$ and $b$ in $\Hom_k(\gr A,k)\cong \Hom_{\gr A}(\gr A\ot \gr A, k)$). The cohomology $\coh^*(A,k)$ is finitely generated over its subalgebra generated by $\xi_a,\xi_b$, as a consequence of the proof of~\cite[Theorem 4.3]{Shroff}. We will see below that the subalgebra generated by $\xi_a,\xi_b$ is in fact a polynomial ring in $\xi_a,\xi_b$, which is Noetherian, so applying the left module version of \cite[Corollary 1.5]{GW}, $\coh^*(A,k)$ is itself (left) Noetherian.

To verify the hypothesis of Theorem \ref{main}, we will want 2-cocycles representing elements in $\coh^*(A\# kG , k)$: We use results in~\cite{NWW}, where the notation is slightly different, with $x$ in place of $a$ and $y$ in place of $b$.
There it is shown directly that there are 2-cocycles $\xi_a,\xi_b$ in $\coh^*(A,k)$
generating a polynomial subring $k[\xi_a,\xi_b]$. Results in \cite[Section~5.1]{NWW} also imply  that $\xi_a,\xi_b$
are in $\ima (\res_{A\# kG,A})$; the needed elements in
$\coh^*(A\# kG , k)$ are constructed explicitly using a twisted tensor product resolution in \cite[Section~3.3]{NWW}. We next claim that $\coh^*(A,k)$ is free with free basis $\{1,\ \eta_a,\ \eta_b,\ \eta_a\eta_b \}$
over the polynomial subalgebra $k[\xi_a,\xi_b]$.\footnote{Since $B$ is a Koszul algebra, $\coh^*(B,k)\cong B^{!}$, the Koszul dual of $B$, which is generated by $\eta_a,\eta_b$ (by abuse of notation) with relations dual to those of $B$, that is, $\eta_a^2= \frac{1}{2}\eta_a\eta_b, \ \eta_b^2=0, \ \eta_b\eta_a = -\eta_a\eta_b$. These relations also hold in $\coh^*(A,k)$, however we do not need this fact.} This will follow once we see that the set 
$$\{\xi_a^i\xi_b^j\eta_a^l\eta_b^m\mid i,j\geq 0, \ l,m = 0,1\}$$
represents  a basis of $\coh^*(A,k)$, since $\xi_a,\xi_b$ commute with each other. Note that the cohomology of $S = \gr A$ is well-known, and has a basis precisely of this form. Recall that Anick's resolution for $A$ is minimal, and a comparison shows that in each degree, the dimensions of $\coh^*(A,k)$ and of $\coh^*(S,k)$ are the same. This forces the May spectral sequence \cite{May} for $A$ to collapse at $E_1 = \coh^*(S,k)$, and so $\gr \coh^*(A,k)\cong\coh^*(S,k)$, and $\coh^*(A,k)$ has basis as claimed. Therefore $\coh^*(A,k)$ is indeed free as a $k[\xi_a,\xi_b]$-module. Further, the $k$-linear span of $\{1,\ \eta_a,\ \eta_b,\ \eta_a\eta_b \}$ is a $kG$-submodule of $\coh^*(A,k)$: We compute
$$ ^g\! \eta_a = \eta_a + \eta_b,  \ \ ^g\! \eta_b = \eta_b, \ \ 
    ^g\!(\eta_a\eta_b) = \eta_a\eta_b.
$$

We have shown that the hypotheses of Theorem~\ref{main} are satisfied. Therefore, $\coh^*(A\# kG,k)$ is Noetherian.

\begin{question}{\em
Are there more examples of Nichols algebras in positive characteristic to which Theorem~\ref{main} applies?
}
\end{question}

Nichols algebras and their bosonizations, which are Hopf algebras, have only just begun to be explored in positive characteristic. There is a vast (and recent) literature on Nichols algebras in characteristic 0. See, for example, \cite{AFGV1,AFGV2,AS,Heckenberger}.


\end{document}